\documentclass[11pt]{article}
\usepackage[francais,english]{babel}
\usepackage[latin1]{inputenc}

\usepackage{amsmath}
\usepackage{amsfonts}
\usepackage{amssymb}
\usepackage{amsthm}
\usepackage{graphics}
\usepackage{amscd}
\usepackage{graphicx,epsfig}
\usepackage{color}

\DeclareMathOperator{\Div}{div}

\renewcommand{\epsilon}{\varepsilon}

\newcommand{\boF}{\mathcal{F}}

\newcommand{\boB}{\mathcal{B}}

\newcommand{\Ome}{\Omega}

\newcommand{\R}{\mathbb{R}}

\renewcommand{\L}{\mathbb{L}}

\newcommand{\N}{\mathbb{N}}

\newcommand{\dd}{\mathrm{d}}

\newtheorem{thm}{Th\'eor\`eme}

\newtheorem{prop}[thm]{Proposition}
\newtheorem{lem}[thm]{Lemme}

\theoremstyle{definition}
\newtheorem{defn}[thm]{D\'efinition}

\renewcommand{\phi}{\varphi}
\newcommand{\dis}{\displaystyle}

\newcounter{remark}

\newcounter{case}

\newcounter{construction}

\newcommand{\arc}[1]{\stackrel{\frown}{#1}}
\DeclareMathOperator{\aire}{Aire}

\title{Lignes de divergence pour les graphes \`a courbure moyenne constante}
\author{Laurent Mazet}
\date{}

\begin{document}
\selectlanguage{french}

\maketitle
\begin{center}
{\small \emph{Universit\'e Paul Sabatier, MIG Laboratoire Emile
  Picard. UMR 5580, 31062 Toulouse cedex 9, France\\
  mazet@picard.ups-tlse.fr}} 
\end{center}

\begin{abstract}
Cet article est consacr\'e \`a l'\'etude de la convergence des suites
de solutions de l'\'equation des surfaces \`a courbure moyenne
constante $H$. Nous d\'efinissons le lieu de convergence de la suite. Le
principale th\'eor\`eme (Th\'eor\`eme \ref{lindiv}) caract\'erise le
compl\'ementaire de ce lieu de convergence: il montre que ce
compl\'ementaire est compos\'e d'arcs de cercle de courbure $2H$. Nous
donnons ensuite des r\'esultats permettant d'exploiter ce th\'eor\`eme.
\end{abstract} 

\selectlanguage{english}
\begin{abstract}
This paper is devoted to the study of convergence of sequences of
solutions to the constant mean curvature $H$ equation. The convergence
domain is defined. The main theorem (Theorem \ref{lindiv})
characterizes the complement of this convergence domain: it shows
that circle arcs of curvature $2H$ compose this complement. We
then give results which allow us to use this theorem.
\end{abstract}
\selectlanguage{french}

\section*{Introduction}
Le but de ce texte est de pr\'esenter une \'etude de la convergence des suites
de solutions de l'\'equation des surfaces \`a courbure moyenne constante:
\begin{equation*}
\label{cmc}
\Div\left(\frac{\nabla u}{\sqrt{1+|\nabla u|^2}}\right)=2H
\tag{CMC}
\end{equation*}

Cette \'equation diff\'erentielle elliptique impose au graphe de $u$ d'\^etre
une surface \`a courbure moyenne constante $H$. Nous nous int\'eressons donc
\`a la convergence des suites $(u_n)$ de solutions de \eqref{cmc} d\'efinies
sur un domaine $\Ome$.

Le cas $H=0$ correspond au cas des surfaces minimales. Pour l'\'equation des
surfaces minimales, une \'etude similaire \`a celle qui va \^etre
pr\'esent\'ee a d\'ej\`a \'et\'e faite dans de pr\'ec\'edents articles
\cite{Ma1,Ma2}.  

Dans la suite, nous allons nous restreindre au cas $H>0$, sachant que le cas
n\'egatif en d\'ecoule si on consid\`ere $-u$ au lieu de $u$. 

Dans un article de 1972 \cite{Sp}, J.~Spruck \'etudie le probl\`eme de
Dirichlet attach\'e \`a l'\'equation \eqref{cmc}. J.~Spruck souhaite autoriser
les valeurs $+\infty$ et $-\infty$ sur le bord du domaine et ainsi trouver un
r\'esultat similaire \`a celui de H.~Jenkins et J.~Serrin pour l'\'equation
des surfaces minimales \cite{JS}. Pour mener ce travail \`a bien, J.~Spruck
s'int\'eresse \`a la convergence des suites monotones de solutions de
\eqref{cmc}. Dans ce cas, la suite converge sur une partie du domaine de
d\'efinition et diverge vers $+\infty$ ou $-\infty$ suivant la
monotonie sur le reste du domaine. J.~Spruck montre que le domaine o\`u la
suite diverge est 
d\'elimit\'e par des arcs de cercle de courbure $2H$. Il obtient gr\^ace \`a
cela des conditions d'existence similaires \`a celles du cas des surfaces
minimales. 

Dans notre article, nous allons montrer que ces arcs de cercle qui d\'elimitent
diff\'erents comportements dans la convergence de la suite existent aussi
lorsque que l'on ote l'hypoth\`ese de monotonie.

Le principale r\'esultat de l'article (th\'eor\`eme \ref{lindiv}) affirme que,
si en un point $P$ du domaine la suite $(\nabla u_n(P))$ n'est pas born\'ee,
il en est alors de m\^eme pour tout point $Q$ appartenant \`a un arc de cercle
de courbure $2H$ de $\Ome$ passant par $P$. Ces arcs de cercles sont appel\'es
lignes de divergence de la suite $(u_n)$.

A l'image de l'\'etude des suites monotones, ce r\'esultat peut \^etre
utilis\'e pour la r\'esolution du probl\`eme de Dirichlet associ\'e \`a
\eqref{cmc}. Il permet aussi d'obtenir des renseignements sur ces solutions
\cite{Ma3}. 

Dans une premi\`ere partie, nous allons pr\'esenter des notations et donner
quelques r\'esultats que nous utiliserons dans les parties suivantes.

Dans la deuxi\`eme partie, nous donnons la d\'efinition du domaine de
convergence d'une suite, il s'agit d'une partie du domaine $\Ome$ o\`u l'on
peut assurer la convergence de $(u_n)$.

La partie \ref{partie3} est consacr\'ee \`a l'\'etude du compl\'ementaire du
domaine de convergence et \`a la d\'emonstration du th\'eor\`eme
\ref{lindiv}. Cette partie aboutit \`a la d\'efinition des lignes de
divergence de la suite $(u_n)$.

La derni\`ere partie contient de nombreux r\'esultats qui permettent
d'exploiter le th\'eor\`eme \ref{lindiv} pour l'\'etude de la convergence ou
de la divergence d'une suite $(u_n)$. Essentiellement, nous \'etudions les
cons\'equences de contraintes sur $\partial\Ome$ par rapport \`a la
convergence de la suite $(u_n)$.


\section{Pr\'eliminaires}
\subsection{Quelques notations}

Dans cette section nous allons fixer quelques notations que nous utiliserons
par la suite.

Consid\'erons $\Ome$ un domaine de $\R^2$ et $P$ un point du bord de
$\Ome$. Comme dans \cite{Se1}, \cite{Se2} et \cite{Sp}, on d\'efinit la
courbure ext\'erieure de $\partial \Ome$ au 
point $P$ comme le supremum des courbures  en $P$ de toutes les courbes $C^2$
passant par $P$ et ne rencontrant pas $\Ome$. Le signe de la courbure est
fix\'e par rapport \`a la normale entrante dans $\Ome$. La courbure
ext\'erieure en $P$ est not\'ee $\hat{\kappa}(P)$. Si aucune courbe $C^2$
n'existe, on pose $\hat{\kappa}(P)=-\infty$. Par exemple, si $\Ome$ est
convexe, en tout point $P\in \partial \Ome$, $\hat{\kappa}(P)\ge 0$.

A plusieurs reprises, nous serons amen\'es \`a consid\'erer des arcs de cercle
de rayon $1/(2H)$. Si $C$ est un tel arc de cercle, sauf mention du contraire,
la normale \`a $C$ que l'on consid\`erera sera toujours celle donn\'ee par le
vecteur de courbure. On notera $\nu$ cette normale, ainsi le vecteur de
courbure sera $2H\nu$ (on rappelle que $H$ est suppos\'e positif). Par
exemple, si $C$ est l'arc de cercle param\'etr\'e par
$$
c:s\longmapsto \left(x_0+\frac{1}{2H}\cos(2Hs), y_0+\frac{1}{2H}\sin(2Hs)
\right) 
$$
avec $s\in I$, la normale est $\nu(c(s))=(-\cos(2Hs),-\sin(2Hs))$.

Par abus de notation, il nous arrivera de voir le vecteur $\nu$ comme un
vecteur horizontal de $\R^3$ en lui adjoignant une troisi\`eme coordonn\'ee
nulle. De m\^eme, il nous arrivera aussi de voir les vecteurs horizontaux de
$\R^3$ comme des vecteurs de $\R^2$. 

\subsection{L'\'equation des surfaces \`a courbure moyenne constante}

Consid\'erons un domaine $\Ome$ de $\R^2$ et $u$ une fonction d\'efinie sur
$\Ome$. Sur le graphe de $u$, la normale que l'on consid\`ere est celle qui
pointe vers le haut. Le graphe de $u$ a alors une courbure moyenne constante
$H$ si $u$ satisfait l'\'equation :

\begin{equation*}
\Div\left(\frac{\nabla u}{\sqrt{1+|\nabla u|^2}}\right)=2H
\tag{CMC}
\end{equation*}


Lorsque l'on \'etudie les solutions de l'\'equation des surfaces \`a courbure
moyenne constante, on introduit la forme diff\'erentielle
$$
\omega_u=\frac{u_x}{W}\dd y -\frac{u_y}{W}\dd x
$$
On utilise les notations suivantes: $u_x$ et $u_y$ d\'esignent les
d\'eriv\'ees premi\`eres de $u$ et $W=\sqrt{1+|\nabla u|^2}$. 
Sur cette forme, l'\'equation \eqref{cmc} se traduit par $\dd
\omega_u= 2H\dd x \wedge\dd y$. On constate ais\'ement que $||\omega_u||\le
1$ ; gr\^ace \`a cela et \`a la remarque pr\'ec\'edente, on peut 
d\'efinir l'int\'egrale de $\omega_u$ sur des arcs inclus dans le bord du
domaine o\`u la fonction $u$ peut ne pas \^etre d\'efinie ou d\'erivable. La
forme diff\'erentielle $\omega_u$ est porteuse de nombreux renseignements sur
la solution $u$ (voir par exemple \cite{Sp}); de plus, on a le lemme suivant.

\begin{lem}\label{lem+inf}
Soit $\Ome$ un domaine du plan et $u$ une solution de \eqref{cmc} sur
$\Ome$. On consid\`ere $C$ un arc du bord de $\Ome$, $C$ est orient\'e comme
faisant partie du bord de $\Ome$ et on d\'esigne par $s$ la longueur d'arc le
long de $C$. On a alors les deux propri\'et\'es suivantes :
\begin{itemize}
\item Si $C$ est un arc de cercle de courbure $\hat{\kappa}=2H$ et
  $\omega_u=\dd s$ le long de $C$, $u$ prend alors la valeur $+\infty$ le long
  de $C$. 
\item Si $C$ est un arc de cercle de courbure $\hat{\kappa}=-2H$ et
  $\omega_u=-\dd s$ le long de $C$, $u$ prend alors la valeur $-\infty$ le
  long de $C$. 
\end{itemize}
\end{lem}

\begin{proof}
Les deux cas se prouvent de fa\c con identique, on ne donne donc que la preuve
du premier point. Consid\'erons $P$ un point de $C$ et notons $D$ l'ensemble
des points de $\Ome$ a distance $r$ de $P$. Pour $r$ suffisament petit, le
r\'esulat de J.~Spruck \cite{Sp} concernant le probl\`eme de Dirichlet pour
\eqref{cmc} nous dit qu'il existe $v$ une solution de \eqref{cmc} sur $D$ qui
prend la valeur $+\infty$ sur la partie du bord de $D$ incluse dans $C$ et les
m\^emes valeurs que $u$ sur le reste du bord de $D$. La
preuve du lemme consiste alors \`a d\'emontrer que $u=v$ sur $D$.

Si $u\neq v$, on peut alors supposer qu'il existe $\eta>0$ tel que
$D'=\{v-u>\eta\}$ soit non-vide. Le bord de $D'$ est compos\'e d'une partie
incluse dans l'int\'erieur de $D$ et \'eventuellement d'une partie incluse
dans $C$. On a $\dd(\omega_v-\omega_u)=0$, donc :
$$
\int_{\partial D'}\omega_v-\omega_u=0
$$
D'apr\`es les hypoth\`eses et le lemme 4.3 de \cite{Sp}, le long de $\partial
D'\cap C$, $\omega_u=\dd s=\omega_v$; ainsi l'int\'egrale sur cette partie du
bord est nulle. Par ailleurs, d'apr\`es le lemme 2 dans \cite{CK},
l'int\'egrale le long du reste du bord de $D'$ est strictement n\'egative car,
sur $D'$, $v\ge u+\eta$. Ceci nous donne une contradiction et donc $u=v$.
\end{proof} 

Remarquons qu'il nous arrivera d'utiliser, par la suite, les notations
usuelles suivantes : $p$ et $q$ pour d\'esigner respectivement $u_x$ et $u_y$
les d\'eriv\'ees premi\`eres de $u$ et $r$, $s$ et $t$ pour $u_{xx}$, $u_{xy}$
et $u_{yy}$. 

\subsection{Une famille de solutions}

Dans cette partie, nous allons pr\'esenter une famille particuli\`ere de
solutions de \eqref{cmc} que nous utiliserons par la suite.

On cherche une solution $u(x,y)$ pr\'esentant une sym\'etrie radiale :
on suppose que $u$ peut s'\'ecrire $u(x,y)=f(r)$ avec
$r=\sqrt{x^2+y^2}$. Si $u$ est solution de l'\'equation des surfaces \`a
courbure moyenne constante, on montre alors que $f$ doit v\'erifier :
$$
\frac{f'}{\sqrt{1+f'^2}}=Hr+\frac{t}{r}
$$
o\`u $t$ est un param\`etre.

Pour que la fonction $f$ existe, il faut et il suffit que $t<1/(4H)$ ; on ne va
s'int\'eresser qu'au cas o\`u $t\in ]0,1/(4H)[$. On montre alors que $f$ n'est
d\'efinie que sur l'intervalle 
$$
r_1(t)=\frac{1-\sqrt{1-4tH}}{2H}\le r\le \frac{1+\sqrt{1-4tH}}{2H}=r_2(t)
$$
et est d\'erivable sur l'intervalle ouvert. On normalise la solution $f$ en
posant $f(1/(2H))=0$. 

Ainsi pour tout $t\in]0,1/(4H)[$, cette fonction $f$ d\'efinit une solution de
\eqref{cmc} que l'on note $h_t$ sur la couronne $r_1(t)\le\sqrt{x^2+y^2}\le
r_2(t)$. Avec la normalisation de $f$, la fonction $h_t$ est nulle sur le
cercle de rayon $1/(2H)$. Le graphe de $h_t$ correspond \`a un morceau
d'ondulo\"\i de d'axe verticale. 

Tout d'abord, on constate que lorsque $t$ tend vers $1/(4H)$, les fonctions
$r_1(t)$ et $r_2(t)$ tendent vers $1/(2H)$.

On constate aussi que la d\'eriv\'ee normale de $h_t$ sur le cercle
$r=r_1(t)$ vaut $-\infty$.

Notons $C$ un arc du cercle d'\'equation $r=1/(2H)$ orient\'e dans le sens
direct. L'expression de $f'$ montre que :
$$
\int_C \omega_{h_t}=(\frac{1}{2}+2Ht)\ell(C) 
$$
o\`u $\ell(C)$ d\'esigne la longueur de l'arc $C$. Ainsi lorsque $t$ tend vers
$1/(4H)$, on obtient :
$$
\lim \int_C \omega_{h_t} = \ell(C)
$$


\section{Le domaine de convergence}
On consid\`ere un domaine $\Ome$ de $\R^2$ et une suite $(u_n)$ de solutions
sur $\Ome$ de l'\'equation \eqref{cmc}. Notre but est d'\'etudier la
convergence \'eventuelle de la suite. L'objectif de cette partie est de
d\'eterminer le lieu o\`u celle-ci converge.

Pour l'\'equation des surfaces \`a courbure moyenne constante, on conna\^it
diff\'erents r\'esultats de convergence, le principale est le suivant. 

\begin{thm}\label{classconv}
Soit $(u_n)$ une suite de solutions de l'\'equation \eqref{cmc} sur un domaine
$\Ome$. On suppose que la suite est uniform\'ement born\'ee sur $\Ome$. Il
existe alors une sous-suite qui converge sur $\Ome$ vers une solution $u$ de
\eqref{cmc}. La convergence est la convergence $C^k$ sur tout compact de
$\Ome$ et ce pour tout $k\in\N$.
\end{thm}

Ce th\'eor\`eme est un r\'esultat classique de compacit\'e pour les solutions
d'\'equations aux d\'eriv\'ees partielles elliptiques. Maintenant, les suites
que nous allons \'etudier ne sont pas, en g\'en\'eral, uniform\'ement born\'ees
sur le domaine ni m\^eme sur tout compact du domaine. Par exemple, si $u$ est
une solution de \eqref{cmc} et $c_n$ est une suite de r\'eels, la suite
$(u+c_n)$ ne converge que si $c_n$ converge. Cette remarque justifie que dans
la suite on s'autorisera des translations verticales pour assurer la
convergence ; une autre fa\c con de voir cela est de dire que l'on
s'int\'eresse essentiellement \`a la convergence de la suite des d\'eriv\'es. 

La convergence qui nous int\'eresse est celle du th\'eor\`eme \ref{classconv}, 
c'est-\`a-dire la convergence $C^k$ sur tout
compact. Entre autre, ceci implique que, si une suite de solutions de
\eqref{cmc} converge, la suites de ses gradients ou la suite $(W_n)$ reste
uniform\'ement born\'ee sur tout compact inclus dans $\Ome$. Cette remarque
est li\'ee \`a notre premier r\'esultat.

\begin{lem}\label{lem}
Soit $\Ome$ un domaine de $\R^2$ et $u$ une solution de l'\'equation des
surfaces \`a courbure moyenne constante sur $\Ome$. On consid\`ere $P$ un
point de $\Ome$ et on note $M= W(P)$. Il existe alors $R>0$ qui ne d\'epend
que de $H$, $M$ et de la distance de $P$ au bord de $\Ome$ tel que, sur le
disque de centre $P$ et de rayon $R$, la fonction $W$ soit major\'ee par $2M$. 
\end{lem}

\begin{proof}
Tout d'abord, on note $r_0$ la distance de $P$ au bord de $\Ome$. La preuve du
lemme repose 
essentiellemnt sur un estim\'e des d\'eriv\'ees secondes de $u$ par 
$W$ qui est d\^u \`a R.~Finn \cite{Fi}. Dans son article, il d\'emontre
(Th\'eor\`eme 2) qu'il existe une constante $C(H,d)$ ne d\'ependant que de $H$
et de la distance $d$ au bord de $\Ome$ telle que, pour tout point $Q$ :
$$
r^2(Q)+s^2(Q)+t^2(Q)\le C(H,d)W^6(Q)
$$ 
L'expression de $C$ est complexe car faisant intervenir des int\'egrales
elliptiques. Toutefois, sa r\'egularit\'e nous permet de dire qu'il existe une
constante $C'$ ne d\'ependant que de $H$ et de $r_0$ telle que, pour tout point
$Q$ dans le disque de centre $P$ et de rayon $r_0/2$, on ait :
$$
r^2(Q)+s^2(Q)+t^2(Q)\le C' W^6(Q)
$$  

Maintenant, on a $\dis\nabla W=(\frac{rp+sq}{W},\frac{sp+tq}{W})$; ainsi sur
le disque $D(P,r_0/2)$ on a $||\nabla W||\le \widetilde{C}W^3$
($\widetilde{C}$ ne 
d\'ependant que de $H$ et $r_0$). Consid\'erons $z$ une solution du probl\`eme
de Cauchy : $z'=\widetilde{C}z^3$ et $z(0)=M$ ; $z$ est d\'efinie sur
$[0,\frac{1}{2M^2\widetilde{C}}[$ par :
$$
\frac{1}{M^2}-\frac{1}{z^2}=2\widetilde{C}r
$$
Ainsi, pour $r\le \frac{3}{8M^2\widetilde{C}}$, $z\le 2M$. Maintenant, la
majoration de $\nabla W$ nous dit qu'en coordonn\'ee polaire autour de $P$:
$W(r,\theta)\le z(r)$. Donc $W$ est major\'e par $2M$ sur le disque
$D(P,\min(r_0/2,\frac{3}{8m^2\widetilde{C}}))$ ; ceci est le r\'esultat
cherch\'e. 
\end{proof}

Nous posons alors la d\'efinition suivante.

\begin{defn}
Soit $\Ome$ un domaine de $\R^2$ et $(u_n)$ une suite de solutions de
\eqref{cmc} sur $\Ome$. On appelle \emph{domaine de convergence} de la
suite $(u_n)$ l'ensemble :
$$
\boB(u_n)=\left\{Q\in\Ome\, |\, (W_n(Q))_{n\in\N} \textrm{ est
    born\'ee}\right\}  
$$
\end{defn}

Le lemme \ref{lem} permet alors de prouver la proposition suivante qui
explique la d\'efinition ci-dessus.

\begin{prop}
Soit $\Ome$ un domaine de $\R^2$ et $(u_n)$ une suite de solutions de
\eqref{cmc} sur $\Ome$. On a  les deux propri\'et\'es suivantes :
\begin{enumerate}
\item L'ensemble $\boB(u_n)$ est un ouvert de $\Ome$.
\item Soit $P\in\boB(u_n)$ et $C$ la composante connexe du domaine de
  convergence contenant $P$, il existe alors une sous-suite de
  $\big(u_n-u_n(P)\big)$ qui converge sur $C$. 
\end{enumerate}
\end{prop}

\begin{proof}
La propri\'et\'e \textit{1} est une cons\'equence directe du lemme
\ref{lem}. Pour la seconde propri\'et\'e, on remarque tout d'abord que le
lemme \ref{lem} implique que la suite $(\nabla u_n)$ est uniform\'ement
born\'ee sur tout compact inclus dans $\boB(u_n)$. Ainsi, sur tout compact
inclus dans $C$, la suite $(u_n-u_n(P))$ est uniform\'ement born\'ee et le
th\'eor\`eme \ref{classconv} permet de construire une sous-suite qui
converge. En consid\'erant une suite exhaustive de compacts de $C$, un
argument de type diagonale de Cantor permet alors de prouver la seconde
propri\'et\'e. 
\end{proof}

La propri\'et\'e \textit{2} illustre l'utilisation de translations verticales
pour assurer la convergence d'une sous-suite. Maintenant la compr\'ehension de
la convergence \'eventuelle d'une suite passe par l'\'etude du
compl\'ementaire du domaine de convergence.


\section{Les lignes de divergence}
\label{partie3}
Le but de cette partie est de comprendre ce qui se passe si un point $P$
n'appartient pas au domaine de convergence d'une suite $(u_n)$ de solutions
de \eqref{cmc}. La premi\`ere \'etape consiste \`a traduire le fait que la
suite $(W_n(P))$ n'est pas born\'ee. Tout d'abord, on remarque que la normale
au graphe au dessus du point $P$ est donn\'ee par :
$$
N_n(P)=\left(\frac{-p_n}{W_n},\frac{-q_n}{W_n},\frac{1}{W_n}\right)(P)
$$
On rappelle que l'on a choisi la normale qui pointe vers le haut. Comme
$(W_n(P))$ n'est pas 
born\'ee, il existe une sous-suite de normales qui converge vers un vecteur
unitaire horizontal. C'est, en fait, cette situation que l'on va essayer de
comprendre. On a un premier r\'esultat.

\begin{prop}\label{conv}
Soit $(u_n)$ une suite de solutions de \eqref{cmc} d\'efinies sur le disque
$D(0,r)$ ; on suppose que la suite des normales au dessus de l'origine
$(N_n(0))$ converge vers le vecteur $(1,0,0)$. On consid\`ere 
$\alpha \in]0,1[$ et on note alors $D_n$ le disque g\'eod\'esique du graphe de
$u_n$ de rayon $\alpha r$  centr\'e en $(0,0,u_n(0))$. Quitte \`a translater
verticalement  $D_n$, une sous-suite de $(D_n)$ converge vers le disque
g\'eod\'esique de rayon $\alpha r$ centr\'e en $(0,0,0)$ du cylindre
vertical d'\'equation $(x-\frac{1}{2H})^2+y^2=(\frac{1}{2H})^2$.
\end{prop}

\begin{proof}
On commence par translater verticalement $D_n$ dans le but de fixer le 
centre des disques en l'origine. Tout point de $D_n$ est \`a distance
sup\'erieure \`a $(1-\alpha)r$ du bord du graphe de $u_n$. Il existe donc une
constante $M$  telle que, pour tout entier $n$ et tout point $Q\in D_n$, la
courbure de Gauss en $Q$ soit born\'ee part $M$ en valeur absolue. La courbure
\'etant born\'ee, il existe une sous-suite de $(D_n)$ qui converge vers $D$ un
disque g\'eod\'esique de centre $(0,0,0)$ et de rayon $\alpha r$ (voir
\cite{MRR}). $D$  est \`a 
courbure moyenne constante $H$ et la normale en son centre est $(1,0,0)$. On
note $N^3$ la troisi\`eme coordon\'ee de la normale \`a $D$. $D$ \'etant la
limite d'une suite de graphe, $N^3$ est positive ; par ailleurs, $N^3$
satisfait l'\'equation suivante : 
$$
\Delta_D N^3=-(4H^2-2K)N^3
$$
Cette \'equation est la traduction de l'harmonicit\'e de l'application de
Gauss sur un surface \`a courbure moyenne constante \cite{He}.
$(4H^2-2K)$ est positif ainsi $N^3$ est sur-harmonique. Or $N^3$ est nulle en
l'origine, elle atteint donc son minimum ; ceci implique que $N^3$ est
constante et donc $N^3=0$. $D$ est donc \`a courbure moyenne constante $H$ et
\`a normale horizontale : ceci implique que $D$ est inclus dans le cylindre
d'\'equation :
$$
(x-\frac{1}{2H})^2+y^2=(\frac{1}{2H})^2
$$
\end{proof}

Cette proposition a pour cons\'equence la convergence de la normale vers des
vecteurs horizontaux au dessus d'un arc de cercle du domaine.

\begin{prop}\label{convpp}
Soit $(u_n)$ une suite de solutions de \eqref{cmc} d\'efinies sur $D(0,r)$ ; on
suppose que $N_n(0)$ converge vers le vecteur $(1,0,0)$. Alors pour tout
$\alpha\in]0,1[$, il existe une sous-suite de $(N_n)$ que l'on note $(N_{n'})$
telle que, pour presque tout $s\in[-\alpha r,\alpha r]$,  on ait :
$$
N_{n'}(c(s)) \longrightarrow (\cos(2Hs),-\sin(2Hs),0)
$$
o\`u $c$ est l'arc de cercle : $\dis c(s)=(\frac{1}{2H}- \frac{1}{2H}
\cos(2Hs), \frac{1}{2H}\sin(2Hs))$.
\end{prop}

\begin{proof}
Consid\'erons $\beta\in]\alpha,1[$. On applique la proposition \ref{conv} aux
disques g\'eod\'esiques $D_n(\beta r)$. Alors, pour une sous-suite que l'on
note $(D_{n'}(\beta r))$, ces disques convergent vers $D$ un disque
g\'eod\'esique inclus dans un 
cylindre vertical de rayon $\frac{1}{2H}$. Ce disque contient, entre autre, la
courbe $\Gamma$ param\'etr\'ee par : $s\mapsto (c(s),0)$ pour $s\in [-\alpha
r,\alpha r]$. Sur $D_{n'}$ et $D$, on d\'efinit les formes diff\'erentielles
$\Ome=\dd X\wedge N$ o\`u $\dd X=(\dd x_1,\dd x_2, \dd x_3)$ et $N$ est la
normale \`a la surface ; on note 
alors $\Ome_3$ la troisi\`eme coordonn\'ee de $\Ome$. Remarquons que ces
formes diff\'erentielles ont une d\'ependence en $n$ que l'on ne marque pas.

Gr\^ace \`a la convergence $D_{n'}(\beta r)\rightarrow D$, il existe
dans le graphe de $u_{n'}$ une courbe $\Gamma_{n'}$ telle que la suite 
$(\Gamma_{n'})$ converge de fa\c con lisse vers $\Gamma$. Plus pr\'ecisement,
on peut, entre autre, assurer que :
\begin{enumerate}
\item si $\gamma_{n'}$ est la projection de $\Gamma_{n'}$ sur le
plan $xy$, alors pour tout $\epsilon>0$, $\gamma_{n'}$ est dans un $\epsilon$
voisinage de $c$ pour $n'$ assez grand ;
\item $\dis\int _{\Gamma_{n'}}\Ome_3\longrightarrow  \int_\Gamma\Ome_3$, ceci
  car la convergence des disques g\'eod\'esiques est r\'eguli\`ere.
\end{enumerate}
On a $\int_\Gamma\Ome_3=-2\alpha r$ ; ainsi $\lim \int
_{\Gamma_{n'}}\Ome_3=-2\alpha r$. Par ailleurs, $\int _{\Gamma_{n'}}\Ome_3 =
\int_{\gamma_{n'}} \omega_{n'}$ o\`u $\omega_n$ est la forme diff\'erentielle
$\omega_{u_n}$ introduite dans la premi\`ere partie. On note $A=c(-\alpha r)$
et $B=c(\alpha r)$ les extr\'emit\'es de $c$ et $A_{n'}$ et $B_{n'}$ les
extr\'emit\'es de $\gamma_{n'}$ correspondantes. Consid\'erons $\epsilon>0$,
pour $n'$ suffisament grand, $\gamma_{n'}$ est dans le
$\epsilon$-voisinage de $c$. On peut alors relier $A$ \`a $A_{n'}$  et $B$ \`a
$B_{n'}$ par des segments de longueur inf\'erieure \`a $\epsilon$. On cr\'ee
ainsi un lacet inclus dans le $\epsilon$-voisinage de $c$ et dont la majeure
partie est compos\'ee de $c$ et 
$\gamma_{n'}$. Comme $\dd \omega_{n'}=2H\dd x\wedge\dd y$, l'int\'egration de
$\omega_{n'}$ le long de ce lacet nous donne la majoration suivante :
$$
\left| \int_c\omega_{n'} - \int_{\gamma_{n'}}\omega_{n'}\right| \le
\left|\int_{[A,A_{n'}]\cup[B,B_{n'}]}\omega_{n'} \right| + 2H
\textrm{Aire}(\textrm{$\epsilon$-voisinage de $c$})=O(\epsilon) 
$$

Cette majoration montre que $\lim \int_c\omega_{n'}=\lim
\int_{\gamma_{n'}}\omega_{n'}= -2\alpha r$. Maintenant on note $N^1$ et $N^2$
les deux premi\`eres composantes de la normale. Alors :
$$
\int_c\omega_{n'}=\int_{-\alpha r}^{\alpha r} -\cos(2Hs)N^1_{n'}(c(s))+
\sin(2Hs)N^2_{n'}(c(s)) \dd s
$$
Or, on a $-\cos(2Hs)N^1_{n'}(c(s))+ \sin(2Hs)N^2_{n'}(c(s))\ge -1$ ; la
convergence de $\int_c\omega_{n'}$ implique donc que
$-\cos(2Hs)N^1_{n'}(c(s))+ \sin(2Hs)N^2_{n'}(c(s))$ converge vers $-1$  dans
$\L^1([-\alpha r,\alpha_r])$. Ceci implique qu'il existe une sous-suite que
l'on indice $n''$ telle que, pour presque tout $s$ :
$$
-\cos(2Hs)N^1_{n''}(c(s))+ \sin(2Hs)N^2_{n''}(c(s))\longrightarrow -1
$$
Ceci se traduit par :
$$
\begin{pmatrix}
N^1_{n''}(c(s))\\
N^2_{n''}(c(s))
\end{pmatrix}
\longrightarrow
\begin{pmatrix}
\cos(2Hs)\\ -\sin(2Hs)
\end{pmatrix}
$$
Ce qui ach\`eve la d\'emonstration.
\end{proof}

Les deux propositions pr\'ec\'edentes sont des r\'esultats locaux ; on en
d\'eduit le r\'esultat global suivant.

\begin{thm}\label{lindiv}
Soit $\Ome$ un domaine et $(u_n)$ une suite de solutions de l'\'equation des
surfaces \`a courbures moyenne constante sur $\Ome$. On consid\`ere $P$ un
point de $\Ome$ et $N$ un vecteur unitaire horizontal. On note $C$ l'arc de
cercle inclus dans $\Ome$ passant par $P$ et de vecteur de courbure $2HN$ en
$P$. Finalement, pour $Q\in C$, on note $\nu(Q)$ le vecteur unitaire horizontal
normal \`a $C$ tel que $2H\nu(Q)$ soit le vecteur de courbure de $C$ en $Q$
($\nu(P)=N$). 

Si la suite de normale $(N_n(P))$ converge vers $N$ alors, pour tout $Q\in C$,
$(N_n(Q))$ converge vers $\nu(Q)$.
\end{thm}

Dans la d\'emonstration, nous allons utiliser la terminologie d'extraction: une
extraction $\theta$ est une application strictement croissante de $\N$ dans
$\N$. Ainsi toute sous-suite de $(u_n)$  peut s'\'ecrire $u_{\theta(n)}$ avec
$\theta$ une extraction. Une sous-extraction $\theta'$ d'une extraction
$\theta$ est une extraction s'\'ecrivant $\theta'=\theta\circ\beta$ avec
$\beta$ une extraction.

\begin{proof}
Tout d'abord on param\'etrise l'arc $C$ par longueur d'arc, on d\'efinit ainsi 
$c:]a,b[\rightarrow \Ome$ avec $c(0)=P$. $C$ est orient\'e de telle fa\c con
que $(\nu(c(s)),c'(s))$ soit une base orthonorm\'ee directe. 

On commence par consid\'erer $\theta$ une extraction, on va alors montrer
qu'il existe une sous-extraction $\theta'$ de celle-ci telle que
$N_{\theta'(n)}(Q)\rightarrow \nu(Q)$ pour presque tout $Q\in C$. Pour cela, on
note $\boF$ l'ensemble des $\epsilon>0$ tel qu'il existe $\theta'$
sous-extraction de $\theta$ avec $N_{\theta'(n)}(Q)\rightarrow \nu(Q)$ pour
presque tout $Q\in c(]a+\epsilon, b-\epsilon[)$. Comme
$N_{\theta(n)}(P)\rightarrow \nu(P)$, la proposition \ref{convpp} 
montre que $\boF$ est non-vide. Soit $\epsilon_0=\inf \boF$, nous allons
montrer que $\epsilon_0=0$ ; pour cela, supposons que $\epsilon_0>0$. On
consid\`ere $P_1=c(a+\epsilon_0)$ et $P_2=c(b-\epsilon_0)$. On choisit $R$ tel
que les disques $D(P_i,R)$ soient inclus dans $\Ome$. Maintenant, d'apr\`es la
d\'efinition de $\epsilon_0$; il existe $Q_1=c(s_1)$ avec $s_1\in
]a+\epsilon_0,a+\epsilon_0+R/3[$, $Q_2=c(s_2)$ avec $s_2\in
]b-\epsilon_0-R/3,b-\epsilon_0[$ et une 
sous-extraction $\theta_1$  de $\theta$ telle que $\lim N_{\theta_1(n)}(Q_i)
=\nu(Q_i)$ et $\lim N_{\theta_1(n)}(Q)=\nu(Q)$ pour presque tout $Q\in
c(]a+\epsilon_0, b-\epsilon_0[)$. On a 
$D(Q_i,2R/3)\in\Ome$, donc on peut appliquer la proposition \ref{convpp} aux
points $Q_i$ avec $\alpha=3/4$. Il existe alors une sous-extraction $\theta_2$
de $\theta_1$ telle que $N_{\theta_2}$ converge vers $\nu$ pour presque tout
point de $c(s_1-R/2,s_1+R/2)$ et $c(s_2-R/2,s_2+R/2)$. Or $s_1-R/2<a
+\epsilon_0$ et $b-\epsilon_0< s_2+R/2$, donc $\epsilon_0$ ne peut \^etre
strictement positif.
$\epsilon_0$ \'etant nul, le proc\'ed\'e diagonal de Cantor permet de
construire la sous-extraction $\theta'$ de $\theta$ souhait\'ee.

Entre autre, il existe des extractions $\theta$ telles que $N_{\theta(n)}$
converge vers $\nu$ pour presque tout $Q\in C$ ; consid\'erons $\theta$ une
telle extraction. On va montrer qu'en fait on a la convergence 
pour tout $Q$ de $C$. Soit $Q\in C$ tel que $N_{\theta(n)}(Q)$ ne converge pas
vers $\nu(Q)$. Tout d'abord, comme dans tout voisinage de $Q$, il existe des
points de $C$ o\`u $N_{\theta(n)}$ converge vers $\nu$, $Q$ n'appartient pas au
domaine de convergence de la suite $(u_n)$. Ainsi il existe $\theta^*$ une
sous-extraction de $\theta$ telle que $N_{\theta^*(n)}(Q)$ converge vers $N^*$
un vecteur unitaire horizontal diff\'erent de $\nu$. A ce vecteur $N^*$ est
associ\'e un arc de cercle $C^*$ de rayon $1/(2H)$ passant par $Q$ comme dans
l'\'enonc\'e du th\'eor\`eme ; on d\'efinit aussi la normale $\nu^*$ le long
$C^*$. D'apr\`es ce que l'on vient de montrer, il existe une
sous-extraction que l'on notera toujours $\theta^*$ telle que
$N_{\theta^*(n)}$ converge vers $\nu^*$ pour presque tout point de $C^*$. 

Quitte \`a changer l'origine de $c$ on peut supposer que $c(0)=Q$. On
param\'etrise $C^*$ par longueur 
d'arc par $c^*$ avec $c^*(0)=Q$ ; $C^*$ est orient\'e comme $C$ par rapport
\`a $N^*$. Soit $\epsilon>0$, on note alors $A=c(-\epsilon)$ et
$B^*=c^*(\epsilon)$. Les arcs de cercles $\arc{AQ}$ et $\arc{QB^*}$ h\'eritent
de l'orientation de $C$ et $C^*$. On consid\`ere alors $D$ le domaine bord\'e
par 
les deux arcs de cercles $\arc{AQ}$ et $\arc{QB^*}$ et le segment $[B^*,A]$,
pour 
$\epsilon$ suffisament petit $D$ est un vrai domaine inclus dans $\Ome$. De
m\^eme avec $B=c(\epsilon)$  et $A^*=c^*(-\epsilon)$, le domaine $D'$ bord\'e
par les arcs de cercles $\arc{A^*Q}$ et $\arc{QB}$ et le segment $[B,A^*]$ est
inclus dans 
$\Ome$ pour $\epsilon$ petit. Suivant les cas, soit l'orientation $AQB^*$
correspond \`a celle de $\partial D$ en tant que bord de $D$ soit c'est celle
de $A^*QB$ qui correspond \`a celle de $\partial D'$ ; on se reporte \`a la
figure \ref{croisement}. 

\begin{figure}[h]
\begin{center}
\resizebox{1\linewidth}{!}{\input{figdiverH2.pstex_t}}
\caption{\label{croisement}}
\end{center}
\end{figure}

On suppose que le premier cas se produit (l'autre cas est identique). On
int\`egre alors $\omega_{\theta^*(n)}$ le long du bord de $D$ :
\begin{equation*}
\begin{split}
\int_{\arc{AQ}}\omega_{\theta^*(n)} + \int_{\arc{QB^*}}\omega_{\theta^*(n)}
+\int_{[B^*,A]}\omega_{\theta^*(n)} &=\int_{\partial D} \omega_{\theta^*(n)}\\
&=2H\textrm{Aire}(D)\\
&> 0
\end{split}
\end{equation*}

Or d'apr\`es la convergence des normales le long de $C$ et $C^*$, on sait que
$\lim \int_{\arc{AQ}}\omega_{\theta^*(n)}=-\epsilon$ et $\lim
\int_{\arc{QB^*}}\omega_{\theta^*(n)}=-\epsilon$. Donc, en utilisant
$||\omega_n||<1$, un passage \`a la limite nous donne $2\epsilon<
\ell([A,B^*])$ ; ceci contredit l'in\'egalit\'e triangulaire.

On sait donc maintenant que $N_{\theta(n)}(Q)$ converge vers $\nu(Q)$ pour tout
$Q\in C$. Soit $Q$ un point de $C$ et supposons que $N_n(Q)$ ne converge pas
vers $\nu(Q)$. Pour une extraction $\alpha$, on peut supposer que
$N_{\alpha(n)}(Q)$ converge vers $N'$ un vecteur unitaire  diff\'erent de
$\nu(Q)$. Or $N_{\alpha(n)}(P)$ converge vers $N=\nu(P)$ donc, d'apr\`es ce
que l'on a  
d\'ej\`a d\'emontr\'e, il existe $\alpha'$ une sous-extraction de $\alpha$
telle que $N_{\alpha'(n)}$ converge vers $\nu$ pour tout point de $C$. Entre
autre, en $Q$, on a $N'=\lim N_{\alpha(n)}(Q)=\lim
N_{\alpha'(n)}(Q)=\nu(Q)$. 

Ceci finit de prouver que $N_n$ converge vers $\nu$ pour tout point de $C$.
\end{proof}

Il y a diff\'erentes cons\'equences que l'on doit retenir de ce
r\'esultat. Tout d'abord, le th\'eor\`eme \ref{lindiv} nous dit que le
compl\'ementaire du domaine de convergence $\boB(u_n)$ est une union d'arc de
cercle de rayon $1/(2H)$. On pose d'ailleurs la d\'efinition suivante.

\begin{defn}
On consid\`ere $C$ un arc de cercle inclus dans $\Ome$ de rayon $1/(2H)$ et on
note $\nu$ la normale \`a $C$ comme dans le th\'eor\`eme \ref{lindiv}.
Soit $(u_n)$ une suite de solutions de \eqref{cmc} d\'efinies sur
$\Ome$. Alors si 
il existe $(N_{n'})$ une sous suite des normales aux graphes des $u_n$
qui converge vers 
$\nu$ pour tout point $Q$ de $C$, on dit que $C$ est une \emph{ligne de
  divergence} de la suite $(u_n)$. On dit de plus que la sous-suite d'indice
$(n')$  fait appara\^itre la ligne de divergence $C$.
\end{defn}

Ainsi le compl\'ementaire du domaine de convergence est l'union des lignes de
divergence de la suite $(u_n)$.

Comme on l'a vu dans les d\'emonstrations, la convergence de la suite des
normales le long d'une ligne de divergence se traduit sur la convergence des
int\'egrales des $1$-formes $\omega_n$. Ainsi, consid\'erons $C$ une ligne de
divergence avec, par exemple, $N_{n_k}(Q)\rightarrow \nu(Q)$ pour tout point
$Q$ 
de $C$. On consid\`ere $T$ un sous-arc de $C$ et on suppose $T$ orient\'e par
un vecteur $v$ tangent en $P\in C$ tel que $(\nu(P),v)$ soit une base
directe. On a alors:
$$
\lim_{k\rightarrow +\infty}\int_T \omega_{n_k}=-\ell(T)
$$
o\`u $\ell(T)$ d\'esigne la longueur de l'arc $T$. D'une mani\`ere
g\'en\'erale, c'est cette caract\'erisation des lignes de divergence qui est
la plus utile.


\section{Quelles lignes de divergence existent ?}
Les deux sections pr\'ec\'edentes nous expliquent les objets que l'on peut
introduire lors de l'\'etude de la convergence d'une suite $(u_n)$ de
solutions de \eqref{cmc}. Toutefois elles ne donnent pas de renseignements qui
permettent de conclure sur la convergence de la suite. Le but de cette section
est de donner des r\'esultat qui permettent cette discussion.

Essentiellement, nous allons donner des r\'esultat qui permettent d'interdire
l'apparition de certaine ligne de divergence. L'id\'ee est qu'une ligne de
divergence a des extr\'emit\'es sur le bord du domaine et donc des conditions
sur $\partial\Ome$ permettent de contr\^oler les lignes de divergence de la
suite. 

\subsection{Le cas des donn\'ees infinies}

Le premier r\'esultat que l'on peut donner concerne le cas o\`u les fonctions
$u_n$ prennent toutes des valeurs infinies le long d'une partie du bord de
$\Ome$. On sait, gr\^ace \`a J.~Spruck, que si une solution $u$ de \eqref{cmc}
prend la valeur $+\infty$ le long d'un arc $A$ du bord de $\Ome$ alors $A$ est
un arc de cercle de courbure ext\'erieure $\hat{\kappa}=2H$. On a alors un
premier r\'esultat.

\begin{prop}\label{+inf}
On consid\`ere $\Ome$ un domaine de $\R^2$ dont un arc $A$ du bord est un arc
de cercle de courbure ext\'erieure $\hat{\kappa}=2H$. On consid\`ere $(u_n)$
une suite de solutions de \eqref{cmc} sur $\Ome$ telle que, pour tout $n$, la
fonction $u_n$ tende vers $+\infty$ en tout point de $A$. Alors aucune ligne de
divergence de la suite $(u_n)$ n'a pour extr\'emit\'e un point int\'erieur \`a
$A$. 
\end{prop}

\begin{proof}
Tout d'abord, supposons qu'une telle ligne de divergence $C$ existe et notons
$P$ l'extr\'emit\'e de $C$ appartenant \`a $A$. On peut supposer que l'arc de
cercle $A$ est l'arc param\'etr\'e par longueur d'arc de la fa\c con suivante
 $a : s\mapsto (\frac{1}{2H}\cos(2Hs), \frac{1}{2H}\sin(2Hs))$ avec $a(0)=P$
 et $-\eta<s<\eta$. On note $\nu$ la normale \`a $C$ telle que, pour tout point
 $Q$ de $C$, $N_{n'}(Q)\rightarrow \nu(Q)$. Quitte \`a sym\'etriser $\Ome$ par
 rapport \`a $y=0$ et consid\'erer la suite $(u_n(x,-y))$ sur le nouveau
 domaine, on peut supposer que $(\nu(P),(-1,0))$ forme une base directe. On
 suivra la suite des notations sur la figure \ref{+infi}.
\begin{figure}[h]
\begin{center}
\resizebox{0.4\linewidth}{!}{\input{figdiverH1.pstex_t}}
\caption{\label{+infi}}
\end{center}
\end{figure}
 On note alors $2\alpha$ l'angle entre la tangente \`a $C$ en $P$ et celle de
 $A$ en 
 $P$, $\alpha$ est inclus dans $]0,\pi/2[$. Consid\'erons $\epsilon>0$, 
 on note $Q_1$ le point de $C$ \`a distance $\epsilon$ de $P$ et $Q_2$ le
 point de $A$ \`a distance $\epsilon$ de $P$ appartenant \`a $y>0$. On
 consid\`ere $D$ le domaine compris entre les arcs de cercle
 $\arc{Q_1P}\subset C$ et $\arc{PQ_2}\subset A$ et le segment $[Q_2,Q_1]$. On
 a alors : 
$$
2H\textrm{Aire}(D)=\int_{\partial D}\omega_{n'}
$$
Ceci implique que :
\begin{equation*}
\begin{split}
\int_{\arc{Q_1P}}\omega_{n'}+\int_{\arc{PQ_2}}\omega_{n'}& =
2H\textrm{Aire}(D) -\int_{[Q_2,Q_1]} \omega_{n'}\\
&\le 2H\textrm{Aire}(D) + \ell([Q_2,Q_1])
\end{split}
\end{equation*}

Commme $u_{n'}$ prend la valeur $+\infty$ le long de $A$,
$\int_{\arc{PQ_2}}\omega_{n'}=\ell(\arc{PQ_2})$ avec $\ell(\arc{PQ_2})$ la
longueur de l'arc $\arc{PQ_2}$. Comme $C$ est une ligne de divergence et que
$N_{n'}(Q)\rightarrow 
\nu(Q)$ pour tout $Q\in C$, un passage \`a la limite dans l'in\'egalit\'e
ci-dessus donne $\ell(\arc{Q_1P})+\ell(\arc{PQ_2})\le 2H\textrm{Aire}(D) +
\ell([Q_2,Q_1])$. Or on sait que
$\ell(\arc{Q_1P})=\ell(\arc{PQ_2})=\epsilon+o(\epsilon)$, $\ell([Q_2,Q_1])\le
2\epsilon\sin\alpha$ et $\textrm{Aire}(D)\le \alpha\epsilon^2$. Donc on a 
$2\epsilon+o(\epsilon)\le 2H\alpha\epsilon^2+2\epsilon\sin\alpha$ ; ceci
implique $2\le 2\sin\alpha$, ce qui est impossible puisque
$\alpha\in]0,\pi/2[$. 
\end{proof}

On a aussi un r\'esultat \'equivalent lorsque toutes les fonctions $u_n$
prennent la valeur $-\infty$ le long du bord. J.~Spruck a montr\'e que, si une
solution $u$ de 
\eqref{cmc} prend la valeur $-\infty$ le long d'un arc $B$ du bord de $\Ome$,
l'arc $B$ est est un arc de cercle de courbure ext\'erieure
$\hat{\kappa}=-2H$. Les techniques de la preuve de la proposition \ref{+inf}
s'adaptent alors pour d\'emontrer le r\'esultat suivant.

\begin{prop}
On consid\`ere $\Ome$ un domaine de $\R^2$ dont un arc $B$ du bord est un arc
de cercle de courbure ext\'erieure $\hat{\kappa}=-2H$. On consid\`ere $(u_n)$
une suite de solutions de \eqref{cmc} sur $\Ome$ telle que, pour tout $n$, la
fonction $u_n$ tende vers $-\infty$ en tout point de $B$. Alors aucune ligne de
divergence de la suite $(u_n)$ n'a pour extr\'emit\'e un point int\'erieur \`a
$B$. 
\end{prop}

Un outil int\'eressant pour \'etudier les lignes de divergence ainsi que les
limites sur le domaine de convergence est donn\'e par le r\'esultat
suivant. Il permet de comprendre le comportement d'une \'eventuelle limite de
la suite $(u_n)$ sur le bord d'une composante du domaine de convergence. 

\begin{prop}\label{valdiv}
On consid\`ere $D(r)$ le disque centr\'e de rayon $r$ et $C$ l'arc de cercle
d'\'equation $(x-1/(2H))^2+y^2=1/(4H^2)$ inclus dans $D(r)$ ($rH$ est suppos\'e
petit). $C$ s\'epare $D(r)$ en deux composantes connexes : l'une contient
$(-r,0)$, elle est not\'e $D^-$ l'autre contient $(r,0)$ elle est not\'ee
$D^+$.

On consid\`ere une suite $(u_n)$ de solutions de \eqref{cmc} d\'efinies
sur $D^-$ qui converge vers une solution $u$ ; on suppose de plus que l'on
satisfait l'une des conditions suivantes :
\begin{enumerate}
\item[1a.] pour tout $n\in \N$, $u_n$ tend vers $-\infty$ sur $C$ ou,
\item[1b.] pour tout $n\in\N$, la fonction $u_n$ est la restriction \`a $D^-$
  d'une  solution $v_n$ de \eqref{cmc} d\'efinie sur $D(r)$ et $C$ est une
  ligne de divergence de $(v_n)$.
\end{enumerate}
Alors la fonction $u$ tend vers $-\infty$ sur $C$.

De m\^eme, si on consid\`ere une suite $(u_n)$ de solutions de \eqref{cmc}
d\'efinies sur $D^+$ qui converge vers une solution $u$ et que l'on suppose de
plus que la suite satisfait l'une des conditions suivantes :
\begin{enumerate}
\item[2a.] pour tout $n\in \N$, $u_n$ tend vers $+\infty$ sur $C$ ou,
\item[2b.] pour tout $n\in\N$, la fonction $u_n$ est la restriction \`a $D^+$
  d'une solution $v_n$ de \eqref{cmc} d\'efinie sur $D(r)$ et $C$ est une
  ligne de divergence de $(v_n)$.
\end{enumerate}
Alors la fonction $u$ tend vers $+\infty$ sur $C$.
\end{prop}

\begin{proof}
Les d\'emonstrations des quatre cas sont semblables, on va donc s'int\'eresser
aux cas \emph{1a} et \emph{1b}. On consid\`ere $C_\epsilon$ l'arc de cercle
inclus dans $D^-$ d'\'equation $(x-(1/(2H)-\epsilon))^2+y^2=1/(4H^2)$. On note 
alors $\Ome_\epsilon$ la partie de $D^-$ comprise entre $C_\epsilon$ et
$C$. On oriente les arcs de $C$ et $C_\epsilon$ dans le sens des $y$
croissant. 

On a alors pour tout $n$ :
$$
\int_{\partial\Ome_\epsilon}\omega_{u_n}=2H \aire({\Ome_\epsilon})
$$
Ceci nous donne, pour tout n :
$$
\left|\int_C\omega_{u_n}-\int_{C_\epsilon}\omega_{u_n}\right|\le
2H\aire(\Ome_{\epsilon})+2 l_\epsilon
$$
o\`u $l_\epsilon$ est la longueur de l'un des deux arcs de cercles qui forment
l'intersection de $\partial \Ome_\epsilon$ et $\partial D(r)$. Dans les deux
cas \emph{1a} et \emph{1b}, $\lim \int_C\omega_{u_n}=-\ell(C)$. Donc  en
passant \`a la limite dans l'in\'egalit\'e ci-dessus, on obtient :
$$
\left|-\ell(C)-\int_{C_\epsilon}\omega_u\right|\le 2H\aire(\Ome_{\epsilon}) +2
l_\epsilon 
$$
Lorsque $\epsilon$ tend vers $0$, $\aire(\Ome_\epsilon)$ et $l_\epsilon$
tendent vers $0$, donc $\int_C\omega_u=-\ell(C)$. Ceci prouve que, le long de
$C$, $\omega_u=-\dd s$. Ainsi, d'apr\`es le lemme \ref{lem+inf}, $u$ prend la
valeur $-\infty$ le long de $C$. 
\end{proof}

\subsection{Le cas des donn\'ees born\'ees}

Dans cette partie, nous allons nous int\'eresser au cas o\`u la suite $(u_n)$
reste finie sur le bord. Dans deux articles \cite{Se1,Se2}, J.~Serrin a
\'etudi\'e le probl\`eme de Dirichlet attach\'e \`a \eqref{cmc} pour des
donn\'ees finies sur le bord du domaine. Il a montr\'e qu'une condition
naturelle pour l'\'etude de ce probl\`eme est de supposer que la courbure
ext\'erieure est partout sup\'erieure \`a $2H$ le long du bord. C'est donc
sous cette hypoth\`ese de courbure que nous allons donner un r\'esultat
concernant les lignes de divergence.

\begin{prop}\label{amelior}
On consid\`ere $\Ome$ un domaine de $\R^2$  dont un arc $C$ du bord a une
courbure ext\'erieure $\hat{\kappa}\ge2H$. On consid\`ere $(u_n)$ une suite de
solutions de \eqref{cmc} sur $\Ome$, continues sur $\Ome\cup C$ telle qu'il
existe $M\in\R$ avec, pour tout $n\in\N$, $|u_n|\le M$ sur $C$. Alors,
aucune ligne de divergence de la suite $(u_n)$ n'a pour extr\'emit\'e un point
int\'erieur \`a $C$.
\end{prop}

\begin{proof}
On remarque tout d'abord que, si $C$ ne contient pas d' arc de cercle de
courbure $\hat{\kappa}=2H$, le lemme 3.3 de \cite{Sp} d\'emontre la
proposition. En effet, celui-ci impose que, dans ce cas, la suite $(u_n)$ soit
born\'ee au voisinage de $C$. Ainsi, d'apr\`es les estim\'es de gradients, la
suite $(W_n(Q))$ reste born\'ee pour tout $Q$ dans un voisinage de $C$ et il
n'y a pas de ligne de divergence. Supposons donc maintenant que $C$ soit un 
arc de cercle de courbure $\hat{\kappa}=2H$. 

D'apr\`es le lemme 3.3 de \cite{Sp}, il existe $c>0$ et un voisinage de $C$
tels que sur ce voisinage $u_n\le M+c$ pour tout $n\in\N$; autrement dit, la
suite est uniform\'ement major\'ee au voisinage de $C$. Supposons que la
suite $(u_n)$ admette une ligne de divergence $C'$ ayant pour extr\'emit\'e
un point $P$ int\'erieur \`a $C$. La suite $(u_n)$ est uniform\'ement
major\'ee sur un voisinage de $P$. Comme dans la d\'emonstration de la
proposition \ref{+inf}, on suppose que l'arc de
cercle $C$ est l'arc param\'etr\'e par longueur d'arc de la fa\c con suivante
$c : s\mapsto (\frac{1}{2H}\cos(2Hs), \frac{1}{2H}\sin(2Hs))$ avec
$c(0)=P$. On note $\nu$ la normale \`a $C'$ telle que, pour tout point $Q$ 
de $C'$, $N_{n_k}(Q)\rightarrow \nu(Q)$; pour simplifier, dans la suite on
supposera que $N_n(Q)\rightarrow \nu(Q)$ . Quitte \`a sym\'etriser $\Ome$ par
rapport \`a $y=0$ et consid\'erer la suite $(u_n(x,-y))$ sur le nouveau
domaine, on peut supposer que $(\nu(P),(-1,0))$ forme une base directe (voir
figure \ref{figuprop}.a). 

La principale difficult\'e de la d\'emonstration est de prouver que, pour $s$
petit, on a :  
\begin{equation*}
\label{butpreuve}
\lim_{n\rightarrow \infty}\int_{c([0,s])}\omega_{u_n}=\ell(c[0,s])=s 
\tag{$*$}
\end{equation*}
On se ram\`ene alors \`a une situation connue.

Pour montrer \eqref{butpreuve}, on va comparer la suite $(u_n)$  \`a des
fonctions barri\`eres que sont les fonctions $h_t$ introduites dans la
premi\`ere partie.  

On va utiliser les coordonn\'ees polaires $(r,\theta)$ de telle fa\c con que
le point $P$ soit le point de coordonn\'ees polaires $(1/(2H), 0)$ et que $C$
soit un arc de $r=1/(2H)$. Soit
$s_0>0$ et $\eta>0$ tels que $c(]-\eta, s_0+\eta[)$ soit inclus dans $C$. L'arc
$C'$ s\'epare un voisinage de $C$ dans $\Ome$ en deux composantes connexes, on
note $\Ome'$ la composante contenant $c(]0,s_0])$ (voir figure
\ref{figuprop}.a).

On consid\`ere $0<t<1/(4H)$ tel que l'ensemble $D=\{(r,\theta)|\,
r_1(t)<r\le1/2H,\ -\eta/2 \le\theta\le s_0\}$ soit inclus dans $\Ome$ et
m\^eme dans un voisinage de $C$ o\`u $u_n\le M+c$ ; pour
$t$ suffisament proche de $1/(4H)$, on a ces inclusions. On va alors comparer
la suite $u_n$ \`a la fonction $h_t-M-1$ qui est d\'efinie sur $D$. On
consid\`ere $\mu_t>0$ tel que $\mu_t$ tende vers $0$ lorsque $t$ tend vers
$1/(4H)$. On pose $s_t$ le plus petit $s\in[\mu_t,s_0]$ tel que l'ensemble
$\{(r,\theta)|\, r_1(t)<r\le1/2H,\ -\eta/2 \le\theta\le s_t/2\}$
contiennent $D\cap C'$. Remarquons que, comme $r_1(t)$ converge vers $1/(2H)$
lorsque $t$ tend vers $1/(4H)$, $s_t$ converge vers $0$ lorsque $t\rightarrow
1/(4H)$. On note alors:
\begin{align*}
D_1&=\{(r,\theta)|\, r_1(t)<r\le1/2H,\ -\eta/2 \le\theta\le s_t\}\\ 
D_2&=\{(r,\theta)|\, r_1(t)<r\le1/2H,\ s_t \le\theta\le s_0\}
\end{align*}
De plus on consid\`ere que $t$ est suffisament proche de $1/(4H)$ pour que
$C'\cap D$ soit un arc de cercle ayant pour extr\'emit\'e $P$ et un point $Q$
contenu dans $\{r=r_1(t)\}$ (voir figure \ref{figuprop}.a). 

\begin{figure}[h]
\begin{center}
\resizebox{0.9\linewidth}{!}{\input{figdiverH3.pstex_t}}
\caption{\label{figuprop}}
\end{center}
\end{figure}

Par hypoth\`ese, $|u_n|<M$ sur $C\cap D$ ; ainsi le long de cet arc, on a
$h_t-M-1=-M-1<u_n$. On va montrer l'\'enonc\'e suivant. 

\begin{lem}\label{chemin}
A partir d'un certain rang $n_0$, il n'existe plus de chemin injectif $\gamma$
dans $D_2$ joignant $C\cap D_2$ \`a $\{r=r_1(t)\}\cap D_2$ le long duquel
$u_n>h_t-M-1$. 
\end{lem}

Supposons au contraire qu'il existe une sous-suite $(n_k)$ d'indice telle que,
pour tout $k$, on puisse trouver un chemin $\gamma_k$ satifaisant la
propri\'et\'e ci-dessus. $\gamma_k$ s\'epare $D$ en deux composantes connexes
dont l'une contient $D_1$, on note $V$ l'intersection de cette composante avec
$\Ome'$ (voir figure \ref{figuprop}.b).

\begin{lem}\label{truc2}
Pour $k$ suffisament grand, il existe des points de $V$ o\`u
$u_{n_k}<h_t-M-1$.
\end{lem}

\begin{proof}
Si ce n'est pas le cas, pour tout $k$, $u_{n_k}\ge h_t-M-1$ sur $D_1\cap
\Ome'$. Ainsi la suite $(u_{n_k})$ est uniform\'ement minor\'ee sur $D_1\cap
\Ome'$ ; on sait par ailleurs que sur cette ensemble elle est major\'ee par
$M+c$. Il existe donc une sous-suite de $(u_{n_k})$ qui converge vers $u$ une
solution de \eqref{cmc} sur $D_1\cap \Ome'$. La fonction $u$ est born\'ee en
tant que limite d'une suite born\'ee. Le domaine $D_1\cap \Ome'$ est en partie
bord\'e par un sous-arc de $C'$ qui est une ligne de divergence ; d'apr\`es le
lemme \ref{valdiv}, $u$ prend la valeur $+\infty$ le long de ce sous-arc. Ceci
contredit le fait que $u$ soit born\'ee et le lemme \ref{truc2} est prouv\'e.
\end{proof} 

On sait maintenant que, pour $k$ grand, il existe des points de $V$ o\`u
$u_{n_k}<h_t-M-1$. Plus pr\'ecis\'ement, on a : 

\begin{lem}\label{truc}
Pour $k$ suffisament grand, il existe des points appartenant
\`a $V\cap \{r=r_1(t)\}$ o\`u $u_{n_k} < h_t-M-1$.
\end{lem}
\begin{proof}
Si le lemme n'est pas v\'erifi\'e, il existe une constante
$K\in\R$ telle que, pour tout $k$, $u_{n_k}$ est minor\'ee par $K$ sur
$\partial V\backslash C'$. Notons $a$ le point de 
coordonn\'ees polaires $(1/(2H),s_0)$ et $b$ celui de coordonn\'ees
$(r_1(t),s_0)$, il s'agit de deux sommets de $D$. Consid\'erons alors le
domaine $U$ de $\R^2$ d\'elimit\'e par les arcs de cercles suivants :
\begin{itemize}
\item le sous-arc $\arc{Pa}$ de $C$,
\item un arc de cercle de courbure $2H$ joignant $a$ \`a $b$
\item un arc de cercle de courbure $2H$ joignant $b$ \`a $Q$
\item un arc de cercle de courbure $2H$ joignant $Q$ \`a $P$
\end{itemize}
Comme sur la figure \ref{figuprop}.c, les arcs sont choisis de telle fa\c con
que par rapport \`a $U$ leur courbure ext\'erieure soit $2H$ pour $\arc{ab}$
et $\arc{bQ}$ et $-2H$ pour $\arc{QP}$ ; on remarque que $U$ contient
$\partial V\backslash C'$. Le 
th\'eor\`eme 6.2 de \cite{Sp} nous dit qu'il existe alors une solution $v$ de
l'\'equation \eqref{cmc}  
d\'efinie sur $U$ tel que $v$ prenne pour valeur sur le bord $0$ le long de
$\arc{Pa}\cup\arc{ab}\cup\arc{bQ}$ et $-\infty$ le long de $\arc{QP}$. La
fonction $v$ est major\'ee par $0$. Ainsi, d'apr\`es le principe du maximum,
pour tout $k\in\N $, $u_{n_k}\ge v+K$ sur $V\cap U \supset D_1\cap U$. Par
ailleurs $(u_{n_k})$ est uniform\'ement major\'ee sur $D_1\cap U$. Ainsi
$D_1\cap U$ est inclus dans le domaine de convergence de la suite
$(u_{n_k})$. 

Regardons maintenant le comportement de la suite sur le domaine de $\R^2$
compris entre le sous-arc de $C'$ joignant $Q$ \`a $P$ et l'arc $\arc{QP}$
bordant $U$. Si un point de ce domaine n'est pas dans le domaine de
convergence, une ligne de divergence doit appara\^itre. Comme le long de $C'$, 
$N_n\rightarrow \nu$, cette ligne de divergence ne peut intersecter $C'$. Ainsi
cette ligne de divergence intersecte l'arc $\arc{QP}$ bordant $U$ et a des
points dans $D_1\cap U$, ce qui est impossible car $D_1\cap U$ est inclus dans
le domaine de convergence. Ainsi on vient de montrer que $D_1\cap \Ome'$ est
inclus dans le domaine de convergence de $(u_{n_k})$. 

Comme $(u_{n_k})$ est uniform\'ement born\'ee sur $D_1\cap U$, quitte \`a
extraire, $(u_{n_k})$ converge vers une solution $u$ de \eqref{cmc} sur
$D_1\cap \Ome'$. Comme $u_{n_k}\le M+c$, $u\le M+c$ sur $D_1\cap \Ome'$. Par
ailleurs le lemme \ref{valdiv} nous dit que $u$ prend la valeur $+\infty$ le
long du sous arc de $C'$ bordant $D_1\cap \Ome'$, ce qui nous donne la
contradiction recherch\'ee et d\'emontre le lemme \ref{truc}.
\end{proof}

\begin{figure}[h]
\begin{center}
\resizebox{0.9\linewidth}{!}{\input{figdiverH4.pstex_t}}
\caption{\label{figuprop2}}
\end{center}
\end{figure}

Reprenons la d\'emonstration du lemme \ref{chemin}. Consid\'erons $k$
suffisament 
grand de fa\c con \`a ce que le lemme \ref{truc} soit satisfait. Le long de
l'arc $D_1\cap C$ on a $u_{n_k}>h_t-M-1$, consid\'erons alors la
composante connexe de $\{u_{n_k}>h_t-M-1\}\cap V$ contenant cet arc et notons
$U$ son compl\'ementaire dans $V$. D'apr\`es les lemmes \ref{truc2} et
\ref{truc}, $U$ est non-vide. 

Comme $u_{n_k}>h_t-M-1$ le long de $\gamma_k$, le bord de $U$ se compose de
trois parties (voir figure \ref{figuprop2}.a) :
\begin{itemize}
\item la premi\`ere, $\Gamma_1$, est l'intersection
de $\partial U$ avec $C'$, il s'agit d'un sous-arc de $C'$ qui est de la forme
$\arc{Q'T}$ si il est non vide,
\item la deuxi\`eme, $\Gamma_2$, est l'intersection avec
$\{r=r_1(t)\}$, il s'agit d'un sous-arc non-vide d'apr\`es le lemme \ref{truc}
et
\item  la derni\`ere, $\Gamma_3$, est la partie incluse dans l'int\'erieur de
  $V$ le long de laquelle $u_{n_k}=h_t-M-1$. 
\end {itemize}
On a alors :
$$
0=\int_{\partial U}\omega_{u_{n_k}}-\omega_{h_t} =\sum_{i=1}^3
\int_{\Gamma_i}\omega_{u_{n_k}}-\omega_{h_t} 
$$

Tout d'abord, le long de $\Gamma_3$ le vecteur $\nabla u_{n_k}-\nabla h_t$
pointe vers l'ext\'erieur de $U$, donc, d'apr\`es le lemme 2 dans \cite{CK},
$\int_{\Gamma_3}\omega_{u_{n_k}}-\omega_{h_t} >0$. Deuxi\`emement, comme la
d\'eriv\'ee normale de $h_t$ en un point de $\{r=r_1(t)\}$ vaut $-\infty$,
$\int_{\Gamma_2}\omega_{u_{n_k}}-\omega_{h_t} >0$. Maintenant $C'$ est une
ligne de divergence et $\Gamma_1$ est de la forme $\arc{QT}$ ainsi pour $k$
suffisament grand $\int_{\Gamma_1}\omega_{u_{n_k}}\ge
\int_{\Gamma_1}\omega_{h_t}$. Donc pour $k$ suffisament grand :
$$
\int_{\partial U}\omega_{u_{n_k}}-\omega_{h_t}>0
$$
Ceci nous donne une contradiction et prouve le lemme \ref{chemin} \qed

Reprenons maintenant la preuve de la proposition \ref{amelior}. On sait que le
long de 
$\{r=1/(2H)\}\cap D_2$ on a $u_n>h_t-M-1$, consid\'erons alors $U$ la
composante connexe de $\{u_n\ge h_t-M-1\}\cap D_2$ contenant cet arc. Le lemme
\ref{chemin} nous dit que pour $n$ suffisament grand le bord de $U$ ne
rencontre pas $\{r=r_1(t)\}$. Ainsi le bord de $U$ se decompose en quatre
parties (voir figure \ref{figuprop2}.b) :
\begin{itemize}
\item l'arc $c([s_t,s_0])=\{r=1/(2H)\}\cap D_2$,
\item un segment $\Gamma_2$ inclus dans $\{\theta=s_0\}$,
\item une courbe $\Gamma_3$ incluse dans $D_2$ joignant $\{\theta=s_0\}$ \`a
  $\{\theta=s_t\}$ le long de laquelle $u_n=h_t-M-1$ et
\item un segment $\Gamma_4$ inclus dans $\{\theta=s_t\}$.
\end{itemize}

Maintenant, on a :
$$
0=\int_{\partial U}\omega_{u_n}-\omega_{h_t}
$$
Cette \'egalit\'e se r\'e\'ecrit :
$$
\int_{c([s_t,s_0])}\omega_{u_n}= \int_{c([s_t,s_0])}\omega_{h_t}-
\int_{\Gamma_3}\omega_{u_n}-\omega_{h_t}-
\sum_{i\in\{2,4\}}\int_{\Gamma_i}\omega_{u_n}-\omega_{h_t} 
$$
Le long de $\Gamma_3$ le vecteur $\nabla u_n-\nabla h_t$ pointe \`a
l'int\'erieur de $U$ donc $\int_{\Gamma_3}\omega_{u_n}-\omega_{h_t}\le 0$. Les
segments $\Gamma_2$ et $\Gamma_4$ sont de longeur inf\'erieur \`a
$1/(2H)-r_1(t)$. On obtient donc l'in\'egalit\'e :
$$
\int_{c([s_t,s_0])}\omega_{u_n}\ge \int_{c([s_t,s_0])}\omega_{h_t}-
4\left(\frac{1}{2H}-r_1(t)\right)  
$$
Cette in\'egalit\'e nous donne :
$$
\ell(c([0,s_0])) \ge \int_{c([0,s_0])}\omega_{u_n}\ge
\int_{c([s_t,s_0])}\omega_{h_t}-
4\left(\frac{1}{2H}-r_1(t)\right)-\ell(c([0,s_t]))
$$ 
Lorsque $t$ tend vers $1/(4H)$, $s_t$ tend vers $0$, $\big(1/(2H)-r_1(t)\big)$
tend vers $0$ et $\int_{c([s_t,s_0])}\omega_{h_t}$ tend vers
$\ell(c([0,s_0]))$. Ainsi en laissant tendre $t$ vers $1/(4H)$, on obtient :
$$
\lim_{n\rightarrow +\infty} \int_{c([0,s_0])}\omega_{u_n}=\ell(c([0,s_0]))
$$
et donc, pour tout $s<s_0$, on obtient \eqref{butpreuve}:
$$
\lim_{n\rightarrow +\infty} \int_{c([0,s])}\omega_{u_n}= \ell(c([0,s]))=s
$$

On se retrouve alors dans une situation similaire \`a celle de
la proposition \ref{+inf} dont on peut appliquer la d\'emonstration avec de
tr\`es l\'eg\`eres modifications. 
\end{proof}

On remarque que la conclusion de la proposition \ref{amelior} est aussi vraie
lorsque l' on a, pour tout $n\in\N$, $\sup_C u_n-\inf_C u_n\le M$ o\`u $C$ est
un arc du bord de courbure ext\'erieure sup\'erieure \`a $2H$. En effet, il
suffit d'appliquer la proposition \ref{amelior} \`a la suite $(u_n-u_n(P))$
o\`u $P$ est un point de $C$.


\end{document}